\documentclass[preprint,12pt]{elsarticle}

\usepackage{amsmath}
\usepackage{amssymb}
\usepackage{dsfont}
\usepackage{amsthm}
\usepackage{mathtools}

\begin{document}
\numberwithin{equation}{section}
\newtheorem{theorem}{Theorem}[section]
\newtheorem{algo}[theorem]{Algorithm}
\newtheorem{corol}[theorem]{Corollary}
\newtheorem{prop}[theorem]{Proposition}
\newtheorem{exempel}[theorem]{Example}
\newtheorem{remark}[theorem]{Remark}
\newtheorem{example}[theorem]{Example}
\newtheorem{conjecture}[theorem]{Conjecture}
\newtheorem{hypo}[theorem]{Hypothesis}
\newtheorem*{theorem*}{Theorem}
\newtheorem{antag}[theorem]{Assumption}
\newtheorem{lemma}[theorem]{Lemma}
\newtheorem{definition}[theorem]{Definition}
\newtheorem{question}[theorem]{Question}
\newtheorem{proposition}[theorem]{Proposition}
\renewcommand{\theequation}{\arabic{section}.\arabic{equation}}
\def\sqr#1#2{{\vcenter{\vbox{\hrule height.#2pt
\hbox{\vrule width.#2pt height#1pt \kern#1pt
\vrule width.#2pt}
\hrule height.#2pt}}}}
\def\ri>{\rightarrow}
\def\->{\longrightarrow}
\def\a{\alpha}
\def\Am{{A_{max}}}
\def\b{\beta}
\def\bR{{\bf R}}
\def\bb{{\cal b}}
\def\cA{{\cal A}}
\def\cB{{\cal B}}
\def\cC{{\cal C}}
\def\cD{{\cal D}}
\def\cE{{\cal E}}
\def\cF{{\cal F}}
\def\cG{{\cal G}}
\def\cL{{\cal L}}
\def\cS{{\cal S}}
\def\cU{{\cal U}}
\def\cX{{\cal X}}
\def\d{\partial}
\def\de{\delta}
\def\e{\varepsilon}
\def\f{\frac}
\def\g{\gamma}
\def\8{\infty}
\def\k{\kappa}
\def\l{\lambda}
\def\L{\Lambda}
\def\N{{\bf N}}
\def\nn{\nonumber}
\def\oa{{\overline a_A}}
\def\oal{{\overline \alpha}}
\def\oB{{\overline B}}
\def\obe{{\overline \beta}}
\def\och{{\overline \chi}}
\def\ocF{{\overline{\cal F}}}
\def\of{{\overline f}}
\def\oF{{\overline F}}
\def\oh{{\overline h}}
\def\ome{\omega}
\def\Ome{\Omega}
\def\ops{{\overline \psi}}
\def\os{{\overline s}}
\def\oS{{\overline S}}
\def\oT{{\overline T}}
\def\ot{{\overline \tau}}
\def\ob{{\overline b}}
\def\obe{{\overline\beta}}
\def\obeta_2{{\overline \beta_2}}
\def\oh{{\overline h}}
\def\ob{{\overline b}}
\def\og{{\overline g}}
\def\oga{{\overline \gamma}}
\def\oi{{\overline i}}
\def\oI{{\overline I}}
\def\oK{{\overline K}}
\def\ok{{\overline k}}
\def\oL{{\overline L}}
\def\om{{\overline m}}
\def\omu{{\overline \mu}}
\def\on{{\overline n}}
\def\oo{{\overline o}}
\def\oq{{\overline q}}
\def\oR{{\overline R}}
\def\oS{{\overline S}}
\def\ote{{\overline t}}
\def\ov{{\overline v}}
\def\ovi{{\overline \vi}}
\def\ow{{\overline w}}
\def\oW{{\overline W}}
\def\ox{{\overline x}}
\def\oX{{\overline X}}
\def\oy{{\overline y}}
\def\oz{{\overline z}}
\def\ophi{{\overline \phi}}
\def\r{\rho}
\def\ovarphi{{\overline \varphi}}
\def\ug{{\underline g}}
\def\uta{{\underline \t}}
\def\0*{^{\odot *}}
\def\qed{\hfill$\sqr45$\bigskip}
\def\rM{{\rm M}}
\def\R{{\mathds R}}
\def\Rn{{\mathds R}^n}
\def\om{\omega}
\def\s{\sigma}
\def\Si{\Sigma}
\def\t{\tau}
\def\th{\theta}
\def\xa{{x_A}}
\def\vi{\varphi}
\def\xb{{x_b}}
\def\xm{{x_m}}
\def\Yx{Y^{\times}}
\def\yx{y^{\times}}
\def\tT{\tilde T}
\def\Proof{\noindent{\em Proof. }}
\def\<{\langle}
\def\>{\rangle}
\begin{frontmatter}


\author[label1]{Istv\'an Bal\'azs}
\ead{balazsi@math.u-szeged.hu}
\author[label2,label4]{Philipp Getto\corref{cor1}}
\ead{phgetto@yahoo.com}
\author[label3,label5]{Gergely R\"ost}
\ead{rost@math.u-szeged.hu}
\fntext[label4]{The research of the author was funded by the
DFG (Deutsche Forschungsgemeinschaft), project number 214819831, and by 
the ERC starting grant EPIDELAY (658, No. 259559).}
\fntext[label5]{The research of the author was funded by the
ERC starting grant EPIDELAY (658, No. 259559), by the Marie Sklodowska-Curie Grant No. 748193, and by NKFIH FK124016.}
 \cortext[cor1]{corresponding author}

\title{A continuous semiflow on a space of Lipschitz functions for a differential equation with state-dependent delay from cell biology}


\address[label1]{Hungarian Academy of Sciences, 1245 Budapest, P.O. Box 1000, Hungary}
 \address[label2]{Center For Dynamics, Technische Universit\"at Dresden, 01062 Dresden, Germany}
\address[label3]{Bolyai Institute, University of Szeged,
Aradi v\'ertan\'uk tere 1, H-6720 Szeged, Hungary; Mathematical Institute, University of Oxford, Woodstock Road, OX2 6GG Oxford, United
Kingdom}


\begin{abstract}
We establish variants of existing results on existence, uniqueness and continuous dependence for a class of delay differential equations (DDE). We apply these to continue the analysis of a differential equation from cell biology with state-dependent delay, implicitly defined as the time when the solution of a nonlinear ODE, that depends on the state of the DDE, reaches a threshold. For this application, previous results are restricted to initial histories
belonging to the so-called solution manifold. We here generalize the results
to a set of nonnegative Lipschitz initial histories which is much larger than the 
solution manifold and moreover convex. 
Additionally, we show that the solutions define a semiflow that  is continuous in the state-component in the $C([-h,0],\R^2)$ topology, which is a variant of established differentiability of the semiflow in $C^1([-h,0],\R^2)$. 
For an associated system we show invariance of convex and compact sets under the semiflow for finite time.
\end{abstract}


\end{frontmatter}
\tableofcontents




\section{Introduction}
 %
%
 %
With this paper we would like to contribute to the development of methods to analyze  differential equations with state-dependent delay  (SD-DDE) 
and to continue the analysis of a model from cell population biology, which can be formulated as a SD-DDE.  In the cell population equation the delay is implicitly defined as the time when the solution of a nonlinear ordinary differential equation meets a threshold (see (\ref{eq11}--\ref{eq14}) below). The SD-DDE additionally features continuously distributed delays. 
%
%

%
%
In \cite{Getto}, the authors have elaborated conditions to
guarantee via application of results of \cite{Walther, Walther1} that the solutions of the 
cell population equation define a differentiable semiflow on the {\it solution manifold}, for $n=2$ a sub-manifold of $\cC^1:=C^1([-h,0],\R^n)$. An advantage of 
the approach in \cite{Walther, Walther1} is the associability of a linear variational equation, from which a characteristic equation, which allows to analyze local stability of equilibria, can be deduced. 
%
%

 %
Motivated by simulations (see the discussion section), a future objective is the proof of existence of periodic solutions for the cell population equation. 
One way to do this is to use fixed point arguments for the
Poincar\'e operator, which is done for a general class of SD-DDE in \cite{MN1}.
As in many fixed point arguments, also in \cite{MN1} convexity and compactness of the domain is used, properties the solution manifold in general does not have. 
Next, note that differentiability of the semiflow in the $\cC^1$-topology as 
established in \cite{Getto} implies continuous dependence on initial values in 
$\cC^1$, i.e., convergence of sequences of solution segments 
in $\cC^1$, if sequences of initial histories converge in $\cC^1$.  The
latter however can appear as too strong in applications, see again the discussion section.
%
%

 %
We here show how - sometimes slightly modified - existing strategies can be combined to show existence, uniqueness and continuous dependence for a large class of SD-DDE. We apply the results to generalize (global) existence and uniqueness of solutions of the SD-DDE (\ref{eq11}--\ref{eq14}) for initial histories in the solution manifold to initial histories in a set of nonnegative Lipschitz functions, the latter being a much larger set than the former and moreover convex. Additionally, we show that the solutions define a semiflow that is continuous in the $\cC:=C([-h,0],\R^n)$ topology. Compared to the above discussed established continuous dependence with respect to initial data in $\cC^1$,  the prerequisite of convergence of initial histories (as well as the conclusion of convergence of solutions) is weaker here  - $\cC$ instead of $\cC^1$ - and we refer to the discussion section for possibilities to exploit this.  
 %
%
%
%

 %
 In \cite{Hale} the existence of noncontinuable and global solutions is established
 for systems of delay differential equations defined by functionals that are continuous on
 domains that are open in the $\cC$-topology ($\cC$-open). Continuous dependence on initial values is shown under the precondition that the solution is unique. Uniqueness
 of solutions is shown if the functional is Lipschitz on a $\cC$-open domain. 
 A known problem is that for SD-DDE the functional is in general not
 Lipschitz on a $\cC$-open domain. A hint to see this is that the evaluation operator
(see (\ref{eq15}) below) is in general not Lipschitz, if functions in the domain are not.
  
 %
 %
%

 %
 %
In \cite{Mallet} the problem is overcome for one-dimensional SD-DDE, where dimension refers to the range space of the functional defining the equation, with the help of the concept of {\it almost local Lipschitzianity}, which roughly means local Lipschitzianity on a domain of Lipschitz functions. It is shown that almost local Lipschitzianity in combination with the discussed results in \cite{Hale} yields existence and uniqueness on a domain of Lipschitz functions.
 %
%
%

 %
 Functionals derived from applications are typically, and in our case, not defined on the whole space but have a domain restricted to a subset of the space.
In  \cite{Mallet} results are first established for an arbitrary functional defined on the whole space $C([-h,0],\R)$.
 Then, to work with restricted domains, a retraction from $C([-h,0],\R)$ to $C([-h,0],[-B,A])$
 is constructed and the results are transferred to the case where the functional
 is defined on $C([-h,0],[-B,A])$ only. A negative feedback condition for the functional
 ensures that solutions remain in the retracted domain. 
 %
 %
%
%

 %
 We here start with a general functional defined on $\cC$. We argue that almost local Lipschitzianity and its use to conclude uniqueness for Lipschitz initial histories can be generalized from one to $n$ dimensions in a straightforward way, conclude uniqueness, and combine it with results from \cite{Hale} on (global) existence and continuous dependence to get existence, uniqueness and continuous dependence for Lipschitz initial histories for a large class of  functionals defined on (all of) $\cC$. 
 %
%
%

 %
To allow for a domain of the form $\cD=C([-h,0],[-B,\infty)^n)$ of the functional, i.e., in particular, a domain that can be specified to our application, we modify the above discussed construction of retractions 
and feedback conditions from \cite{Mallet}. One then can work with a retraction from $\cC$ to $\cD$ and a component-wise feedback condition and transfer the general results on solutions to the case of a functional defined on $\cD$. We conclude that the solutions define a semiflow, in the sense of e.g. \cite{Amann}, that is continuous in the $\cC$-topology on a set of Lipschitz functions and use this continuity to derive some
further properties. 
%
%
%

 %
 We then establish compactness results employing the following ideas. In \cite{MN1} it is used that by the Arzela-Ascoli theorem a set of functions that share the same finite bound and finite Lipschitz constant is compact in $\cC$.  As will be motivated, the approach of \cite{MN1} to show that a time $t$ map leaves such a set invariant for arbitrarily
large  $t$ does not work here directly. However, a class of two-dimensional systems that contains (\ref{eq11}--\ref{eq14}) can be transformed
to a one-dimensional equation. For the latter, invariance of a compact set for finite time can be elaborated. We refer to the discussion section for more details on future implementation of these results.   
%
%

 %
 %
After having established the general results, we consider the SD-DDE
\begin{eqnarray}
w'(t)&=&q(v(t))w(t),
\label{eq11}\\
v'(t)&=&
\frac{\g(v(t-\t(v_t)))}{g(x_1,v(t-\t(v_t)))}g(x_2,v(t))w(t-\t(v_t))e^{\int_0^{\t(v_t)}
[d-D_1g](y(s,v_t),v(t-s))ds}
\nn\\
&&-\mu v(t),
\label{eq12}
\end{eqnarray}
where $y=y(\cdot,\psi)$ and $\t=\t(\psi)$ are defined as the respective solutions of
\begin{eqnarray}
&&y'(s)=g(y(s),\psi(-s)),\;s>0,\;\;y(0)=x_2\;\;{\rm and}
\label{eq13}\\
&&y(\t,\psi)=x_1, 
\label{eq14}
\end{eqnarray}
where $x_1<x_2$ are given parameters. As common in delay differential equations (DDE) we use the notation $x_t(s):=x(t+s)$, $s<0$, for functions $x$ defined in
$t+s\in\R$.
The system describes the dynamics of a stem cell population ($w$) regulated by the mature cell population ($v$). We refer to \cite{Getto} and references therein, in particular \cite{Getto1}, for biological background of the model. 
The SD-DDE can be deduced via integration along the characteristics from a partial
differential equation of transport type which features a progenitor cell maturity density
and maturity structure, see \cite{Getto}. 
%
%
We apply our general results to (\ref{eq11}--\ref{eq14}). 
%
%
To guarantee some of the required conditions, we show that the functional that defines the 
system is almost locally Lipschitz. To handle the implicitly defined state-dependent delay 
we consider evaluation operators and implicitly defined operators and analyze them on 
Lipschitz subsets of continuous functions. 

The paper is structured top down: In Section \ref{s2} we consider our most general class of equations. Section \ref{s3} contains results for an intermediate class and Section \ref{s1} an application of these results to the stem cell SD-DDE; in each of these two sections a subsection on main results precedes one on proofs. Finally, Section \ref{ss5} contains examples of modelling ingredients and Section \ref{s6} a discussion of our results and potential future applications.

%
%
\section{Solving DDE on a state space of Lipschitz functions}\label{s2}
%
%
 \subsection{Initial value problem}\label{ss0}
\begin{definition}\rm
Suppose that $\phi\in\cD\subset\cC$ and $f:\cD\->\Rn$. By  a {\it solution} of
\begin{eqnarray}
x'(t)&=&f(x_t),\;\;\;t\ge t_0,
\label{eq7}
\\
x_{t_0}&=&\phi,
\label{eq1}
\end{eqnarray}
or a solution of (\ref{eq7}) through $\phi$,
we mean a continuous function $x^\phi:[t_0-h,t_0+\a]\->\Rn$ for some $\a>0$,
that is such that on $[t_0,t_0+\a]$ one has $x_t^\phi\in\cD$, the function $x^\phi$ is differentiable and satisfies (\ref{eq7}--\ref{eq1}).
 Solutions on half-open intervals 
$[t_0-h,t_0+\a)$ for $\a\in(0,\infty]$ are defined analogously. 
\end{definition}
We shall sometimes write $x$ instead of $x^\phi$.
\subsection{Domain of the functional is $\cC$}\label{ss3}
%
 %
 %
%
%
\subsubsection{Noncontinuable and global solutions}
\begin{theorem}\label{theo3} 
Suppose that $F: \cC\->\Rn$ is continuous and $\phi\in \cC$. Then
\begin{itemize}
\item[(a)] there exists a unique $c=c(\phi)\in(0,\infty]$ such that $x^\phi:[t_0-h,t_0+c)\->\Rn$ is a non-continuable solution of
\begin{eqnarray}
x'(t)=F(x_t),\;t\ge t_0,\;\;x_{t_0}=\phi.
\label{eq6}
\end{eqnarray}
\end{itemize}
If additionally $F(U)$ is bounded whenever $U\subset \cC$ is closed and bounded then the following hold: 
\begin{itemize}
\item[(b)] If $c<\infty$ then for any closed and bounded $U\subset \cC$ there exists some $t_U\in(0,c)$
such that $x_t^\phi\notin U$ for all $t\in [t_0+t_U, t_0+c)$. 
\item[(c)] If $\{x_t^\phi: t\in[t_0,t_0+\a)\}\subset \cC$ is bounded, whenever $\a<\infty$ and $x^\phi$ is defined on
$[t_0,t_0+\a)$, then $c=\infty$, i.e., the solution is global. 
\end{itemize}
\end{theorem}
The existence of a solution $x^\phi:[t_0-h,t_0+\a]\->\Rn$ for some $\a>0$ 
follows from \cite[Theorem 2.2.1]{Hale} and the statement in (a) is concluded in \cite[Section 2.3]{Hale} from Zorn's lemma. Next, (b) follows from \cite[Theorem 2.3.2]{Hale}. Then (c) is standard: If $c<\infty$ define 
\[
U:=\overline{\{x_t^\phi:\;t\in[t_0,t_0+c)\}}.
\] 
Then by (b) there exists some $t_U\in(0,c)$ such that $x_{t_0+t_U}^\phi\notin U$, which contradicts the definition of $U$. 
\begin{remark}
Note that the cited results in \cite{Hale} hold for non-autonomous equations. Since
our motivation here is an autonomous system and moreover the uniqueness result that we
will use is also for autonomous systems we have rewritten these results
for the autonomous case. 
\end{remark}
%
%
\subsubsection{Uniqueness}

To guarantee uniqueness, the notion of almost local Lipschitzianity for $n=1$ from \cite{Mallet} can
be generalized to arbitrary finite dimensions in a straightforward way.
As common, we define for any $\phi\in\cC$
\[
lip\;\phi:=\sup\left\{\frac{|\phi(s)-\phi(t)|}{|s-t|}:\;s,t\in[-h,0],\;s\neq t\right\}\in[0,\infty]
\]
and $B_\de(x_0):=\{x:\;\|x-x_0\|<\de\}$, where $\de>0$, $|\cdot|$ denotes norms in $\R^n$
with $n$ depending on context,  and the choice of norm $\|\cdot\|$ should
also be clear from the context, e.g., the choice of $x_0$. In the following, however, we denote by $\|\cdot\|$ the sup-norm on $\cC$. Then, a function $\phi$ is Lipschitz
with Lipschitz constant $k$ (we will write $k$-{\it Lipschitz}) whenever $\infty>k\ge lip\;\phi$. 
For each $\phi_0\in \cC$, $\de>0$, $R>0$ define
\[
V(\phi_0;\de,R):=\{\phi\in\oB_\de(\phi_0):\;lip\;\phi\le R\}.
\]
\begin{definition}\label{def1}\rm
A functional $f:\cD\subset \cC=C([-h,0],\R^n)\->\R^m$ is called {\it almost locally Lipschitz} if $f$ is continuous and 
for all $\phi_0\in\cD$, $R>0$ there exists some
$\de=\de(\phi_0,R)>0$, $k=k(\phi_0,R,\de)\ge 0$ such that for all $\vi,\psi\in V(\phi_0;\de,R)\cap\cD$
\[
|f(\vi)-f(\psi)|\le k\|\vi-\psi\|. 
\]
\end{definition}
The following theorem is proven as \cite[Theorem 1.2]{Mallet} for the case $n=1$. The proof
for general $n$ is analogous and we omit it. 
For $\cD\subset \cC$, define $V_\cD:=\{\phi\in\cD:\;lip\;\phi<\infty\}$. Note that if $\cD$ is
convex, so is $V_\cD$. 
\begin{theorem}\label{theo2}
Suppose that $F:\cC\->\Rn$ is almost locally Lipschitz. Let $\phi\in V_\cC$ and $t_0\in\R$. If $\a>0$ and $y,z:[t_0-h,t_0+\a]\->\Rn$ are both solutions of  (\ref{eq6}),
then $y(t)=z(t)$ for all $t\in[t_0,t_0+\a]$. 
\end{theorem}
%
%
\subsubsection{Continuous dependence on initial values}
The following result follows directly from \cite[Theorem 2.2.2]{Hale} if we use our uniqueness 
result. 
\begin{theorem}\label{theo5}
Suppose that $F:\cC\->\R^n$ is almost locally Lipschitz, $\phi\in V_\cC$
and let $\a>0$ be such that  a solution $x^\phi$ through $\phi$ exists on $[t_0-h,t_0+\a]$. Let $(\phi^k)\in V_\cC^\N$ with
$\phi^k\->\phi$. Then $x^\phi$ is unique on $[t_0-h,t_0+\a]$, for some $k\ge k_0$
there exist unique solutions $x^k$ through $\phi^k$ on $[t_0-h,t_0+\a]$ for all $k\ge k_0$ and $x^k\->x^\phi$ uniformly on $[t_0-h,t_0+\a]$. 
\end{theorem}
\begin{remark}
Note that similarly as in \cite[Theorem 2.2.2]{Hale} we could include continuous dependence
on functional and initial time in the above formulation. We did not do this, since, especially
when transferring these results to restricted domains, the exposition would suffer from 
further technicalities and moreover we currently see no direct use for these properties. 
\end{remark}
%
%
%
\subsection{Retraction onto a specific domain}\label{ss4}
It is remarked in \cite{Mallet} (without proof) that the following result holds in case $n=1$. The proof for general $n$ is analogous and we present it for completeness. 
Recall that a {\it retraction} is a continuous
map of a topological space into a subset that on the subset equals the identity. 
\begin{lemma}\label{lem1}
Let $\cD\subset \cC$, $\r :\cC\->\cD$ be a locally Lipschitz retraction. Suppose that for all $\phi_0\in \cC$, $\de>0$, $R>0$
\[
\sup\{lip\;\r(\phi):\;\phi\in V(\phi_0;\de,R)\}<\infty.
\]
Then, if $f:\cD\->\R^n$ is almost locally Lipschitz, so is $F:\cC\->\R^n;F:=f\circ\r$. 
\end{lemma}
\Proof First, $F$ is continuous as a composition of continuous functions. Next, let $\phi_0\in\cC$, $R>0$. Define $L:=\sup\{lip\;\r(\phi):\;\phi\in V(\phi_0;1,R)\}<\infty$. Choose $\e=\e(\r(\phi_0),L)$, $k=k(\r(\phi_0),L)$ such that $f$ is $k$-Lipschitz on $V(\r(\phi_0);\e,L)$. Choose $\de<1$, 
$l\ge 0$ such that $\r(B_\de(\phi_0))\subset B_\e(\r(\phi_0))$ and $\r$ is $l$-Lipschitz
on $B_\de(\phi_0)$. Then for $\vi,\psi\in V(\phi_0;\de,R)$, one has
\begin{eqnarray}
|F(\vi)-F(\psi)|=|f(\r(\vi))-f(\r(\psi))|\le k|\r(\vi)-\r(\psi)|\le kl\|\vi-\psi\|. 
\nn
\end{eqnarray}
Hence, $F$ is $kl$-Lipschitz on $V(\phi_0,\de,R)$ and thus almost locally Lipschitz. 
\qed

%
\subsubsection{A specific retraction for a specific domain}
For the remainder of the section we will use the following construction (unless specified
otherwise). 
\begin{remark}
The construction is a modification of the retraction in \cite{Mallet}, the latter of which maps $C([-h,0],\R)$ onto $C([-h,0],[-B,A])$ with $-\infty<-B<A<\infty$, to a retraction
of $C([-h,0],\R^n)$ onto $C([-h,0],[-B,\infty)^n)$ with $-\infty<-B$. With the result we can work with
nonnegative solutions, if $B=0$, of multi-dimensional systems. The construction could probably be generalized to the range 
$C([-h,0],\Pi_{i=1}^{n}[-B_i,A_i])$, $-\infty\le -B_i <A_i\le\infty$, $ i=1,...,n$.
\end{remark}
Let $B\in\R$ and define
\begin{eqnarray}
\cD:=C([-h,0],[-B,\infty)^n).
\label{eq2}
\end{eqnarray}
Note that the convexity of $\cD$ implies convexity of $V_\cD$. 
We define a map 
\begin{eqnarray}
r:\R\->[-B,\infty),\;
r(u):=
\begin{cases}
u,& u\in[-B,\infty),
\\
-B,& u<-B. 
\end{cases}
\label{eq3}
\end{eqnarray}
Then $r$ is a retraction and Lipschitz with $lip\;r\le 1$. With $r$ we define
another map 
\begin{eqnarray}
\r:\cC\->\cD, \r=(\r_1,...,\r_n), 
\r_i(\phi)(t):=r(\phi_i(t)),\;
i=1,...,n. 
\label{eq4}
\end{eqnarray}
\begin{lemma}
 $\r$ is a retraction and maps bounded sets into bounded sets. 
 \end{lemma}
\Proof It is clear that $\r$ (is onto,) preserves the subset and maps bounded sets into bounded sets. Regarding continuity, suppose
that $\phi^n\->\phi$, and let $\e>0$. Then
\begin{eqnarray}
|[\r_i(\phi^n)-\r_i(\phi)](t)|=|r(\phi^n_i(t))-r(\phi_i(t))|.
\nn
\end{eqnarray}
Choose $N\in\N$, $\de>0$ such that $\|\phi^n-\phi\|\le\de$ for all $n\ge N$. Then
\[
|\phi^n(t)|\le\|\phi\|+\de,\;|\phi(t)|\le\|\phi\|+\de,\;\forall t\in[-h,0], \;n\ge N. 
\]
Now, continuity follows by uniform continuity of $r$ on compact sets. 
\qed

The following result follows by definition of $\r$ from Lipschitzianity of $r$ with $lip\;r\le 1$. We omit the straightforward proofs of the two following results. 
\begin{lemma}\label{lem15}
One has $lip\;\r(\phi)\le lip\;\phi$, hence if $\phi$ is Lipschitz so is $\r(\phi)$. Moreover, $\r$ is Lipschitz with $lip\;\r\le 1$. 
\end{lemma}
The result implies that $\sup\{lip\;\r(\phi):\;\phi\in V(\phi_0;\de,R)\}\le R<\infty$ for all $\phi_0\in \cC$, $\de>0$, $R>0$.
We can use the latter to directly apply Lemma \ref{lem1} to $F:=f\circ\r$ with $\r$ 
being our (locally) Lipschitz retraction:
 \begin{lemma}\label{lem2}
Suppose that $f:\cD\subset\cC\->\Rn$ is almost locally Lipschitz. Then so is $F$. 
 \end{lemma}
 %
 %
 \subsubsection{Noncontinuable and global solutions and uniqueness}
To guarantee that a solution remains within a domain a feedback condition can be used. The proof of the following result is a modification of a 
similar result for one dimension \cite[Theorem 1.3]{Mallet}. 
\begin{lemma}\label{lem17}
Suppose that $f:\cD\->\R^n$ satisfies
\[
f_i(\phi)\ge 0,\;{\rm if}\;\phi_i(0)=-B,\;\forall\phi=(\phi_1,...,\phi_n)\in\cD,\;i=1,...,n
\tag{F}.
\]
Now fix $\phi\in\cD$ and assume that $x$ is a solution of $x'(t)=f(\r(x_t))$ through $\phi$ on some
interval $[t_0-h,t_0+\a]$. Then $x_t\in\cD$  and thus $\r(x_t)=x_t$ for all $t\in[t_0,t_0+\a]$ and hence $x$ is a solution of (\ref{eq7}--\ref{eq1}) on $[t_0,t_0+\a]$. 
\end{lemma}
\Proof The statement would follow if 
$x_i(t+\th)\ge-B$ for all $t\ge t_0$, $\th\in[-h,0]$, $i=1,...,n$. First, $\phi\in\cD$ implies
that $\phi_i(\th)\ge-B$ for all $\th\in[-h,0]$, $i=1,...,n$. Suppose that for some $i\in\{1,...,n\}$
and $x=x^\phi$ one has $x_i(t_1)<-B$ for some $t_1>t_0$. Then $\t:=\sup\{t\in[t_0,t_1]:
\;x_i(t)=-B\}\in[t_0,t_1)$. Then $x_i(\t)=-B$, $x_i(t)<-B$ for all $t\in(\t,t_1]$. By the mean value
theorem $x_i'(t)<0$ for some $t\in(\t,t_1)$. Then $\r_i(x_t)(0)=r(x_i(t))=-B$. Hence by (F)
we have $x_i'(t)=f_i(\r(x_t))\ge 0$, which is a contradiction. 
\qed

 \begin{theorem}\label{theo6} 
Suppose that $f:\cD\->\Rn$ is continuous and satisfies (F). Then the following hold. 
\begin{itemize}
\item[(a)] For every $\phi\in\cD$ there exists a unique $c=c(\phi)\in(0,\infty]$ and a non-continuable solution
$x^\phi$ on $[t_0-h,t_0+c)$ of (\ref{eq7}--\ref{eq1}). 
\end{itemize}
\begin{itemize}
\item[(b)] If $f(U)$ is bounded, whenever $U\subset \cD$ is bounded, and if for some
$\phi\in \cD$ the set $\{x_t^\phi:\;t\in[t_0,t_0+\a)\}\subset \cD$ is bounded, whenever
$\a<\infty$ and $x^\phi$ defined on $[t_0,t_0+\a)$, then $c=\infty$, i.e., the solution is global. 
\item[(c)] If $f$ is almost locally Lipschitz and $\phi\in V_\cD$, then $x^\phi$ is unique.
\end{itemize}
 \end{theorem}
\Proof  
Since $F:=f\circ\r$ is continuous, by Theorem \ref{theo3} (a) there exists a noncontinuable
solution of (\ref{eq6}) for this $F$.  Next, suppose
that $U\subset \cC$ is (closed and) bounded. Then, as remarked, $\r(U)\subset\cD$ is bounded
and hence by the assumption of (b) $F(U)=f(\r(U))$ is bounded. Thus by Theorem \ref{theo3} (c) we have shown
that if $\{x_t^\phi:t\in[t_0,t_0+\a)\}\subset \cC$ is bounded
whenever $\a<\infty$ and $x^\phi$ defined on $[t_0,t_0+\a)$, then $c=\infty$. 
If $f$ is almost locally Lipschitz, then by Lemma \ref{lem2} so is $F$ and thus 
by Theorem \ref{theo2} we get uniqueness. To complete the proof note that (F) guarantees via Lemma \ref{lem17} that $\{x_t^\phi:t\in[t_0,t_0+\a)\}\subset \cD$ and that 
$x^\phi$ is a solution of (\ref{eq7}--\ref{eq1}). 
\qed

\begin{remark}
If $f$ would map only the closed and bounded sets on bounded sets, as required in Theorem
\ref{theo3}, we could not guarantee that $F(U)=(f\circ\r)(U)$ is bounded if $U$ is 
closed and bounded: the above defined retraction $\r$ maps bounded on bounded, 
but in general does not map closed and bounded on closed sets. To see the latter, 
consider e.g. $\cC:=C([0,2],\R)$, $\cD:=\{x\in\cC:x(t)\ge0,\;\forall\;t\in[-h,0]\}$ and $r$ and $\r$ defined as above, but for $n=1$, $B=0$ and the modified $\cC$ and $\cD$. 
Define $U:=\{x_n:\;n\ge 2\}\subset\cC$, where 
\begin{eqnarray}
x_n(t):=\begin{cases}
\frac{1}{n}&,t<1-\frac{1}{n}
\\
1-t&,1-\frac{1}{n}\le t<1
\\
-n(t-1)&,1\le t<1+\frac{1}{n}
\\
-1&,1+\frac{1}{n}\le t\le 2. 
\end{cases}
\nn
\end{eqnarray}
Then it is easy to see that $U$ is closed and bounded but 
\[
\r(U)=\{x:\;\exists\;n\ge2\;s.th.\;x(t)=x_n(t)\;\forall\;t\in[0,1],\;x(t)=0\;\forall\;t\in[1,2]\}. 
\]
is not closed. 
\end{remark}
%
%
\subsubsection{Continuous dependence on initial values}
The negative feedback condition $(F)$ now ensures that our results on continuous 
dependence can be transferred to our case of a specific retraction onto the domain of
the functional. 
\begin{theorem}\label{theo7}
Suppose that $f:\cD\->\R^n$ is almost locally Lipschitz and satisfies (F), let $\phi\in V_\cD$ and $\a>0$ be such that a solution $x^\phi$ of (\ref{eq7}--\ref{eq1}) through $\phi$ exists on $[t_0-h,t_0+\a]$.
Let $(\phi^k)\in V_\cD^\N$ with $\phi^k\->\phi$. Then $x^\phi$ is unique on $[t_0-h,t_0+\a]$, for some $k\ge k_0$ there exist unique solutions $x^k$ through 
$\phi^k$ on $[t_0-h,t_0+\a]$ and
$x^k\->x^\phi$ uniformly. 
\end{theorem}
\Proof Since $x^\phi$ is a solution of (\ref{eq7}--\ref{eq1}) we have $x_t^\phi\in\cD$ for all
$t\ge t_0$. Thus, for $F:=f\circ\r$, one has $F(x^\phi_t)=f(x^\phi_t)$ and $x^\phi$ is  a solution of $x'(t)=F(x_t)$ through $\phi$. Since
$f$ is almost locally Lipschitz, by Lemma \ref{lem2} so is $F$ and since $\phi\in V_\cD$ the solution is unique. By Theorem \ref{theo5} there exists some $k_0$, such that for all $k\ge k_0$ there exist unique solutions
$x^k$ of $x'(t)=F(x_t)$ through $\phi^k$ on $[t_0-h,t_0+\a]$ and $x^k\->x^\phi$ uniformly. By Lemma
\ref{lem17} we have $x_t^k\in\cD$ for all $t\ge t_0$, hence the $x^k$ solve also (\ref{eq7}--
\ref{eq1}). 
\qed

\subsubsection{A continuous semiflow on a state-space of Lipschitz functions}
If $f$ satisfies the assumptions for global existence and uniqueness, we can use the concept of a semiflow, e.g., in the sense of \cite[Section 10]{Amann}.
We start with some definitions:
\begin{definition}\rm
Let $(X,d)$ be a metric space. A map $\Si:[0,\infty)\times X\-> X$ is called a {\it continuous
semiflow} if
\begin{itemize}
\item[(i)] $\Si(0,x)=x$ for all $x\in X$,
\item[(ii)] $\Si(t,\Si(s,x))=\Si(t+s,x)$ for all $s,t\in[0,\infty)$, $x\in X$ (``semigroup property''),
\item[(iii)] $\Si$ is continuous.
\end{itemize}
A {\it trajectory} of the semiflow $\Si$ is a map $\s:I\-> X$ defined on an interval $I\subset \R$
with positive length, such that for $s$ and $t$ in $I$ with $s\le t$ one has 
\[
\s(t)=\Si(t-s,\s(s)). 
\]
The $\ome$-limit set of a trajectory $\s:I\-> X$ with $\sup I=\infty$ is defined as
\[
\ome(\s)=\{x\in X:\;\exists\;(t_n)\in I^\N,\;{\rm s.th.}\;t_n\->\infty,\;\s(t_n)\-> x\;{\rm as}\;n\ri>\infty\}. 
\]
\end{definition}
\begin{remark}
Note that the definitions in \cite[Section 10]{Amann} and \cite[Definition VII 2.1]{3Diekmann} include also semiflows induced by local solutions. Moreover  \cite[Definition VII 2.1]{3Diekmann} additionally requires completeness of the metric space,
which we here cannot expect, since by the Weierstrass approximation theorem 
$V_\cD$ is not complete. On the other hand to our understanding this 
completeness is not  necessary here. 
Note also that \cite[Definition VII 2.1]{3Diekmann} merely requires continuity in each of the components, point-wise with respect to the other. 
The definitions of trajectories
and $\ome$-limit sets are consistent with \cite[Definitions VII 2.3 and 2.4]{3Diekmann}. Note that the reference also contains similar results for $\a$-limit sets. 
\end{remark}
The following properties of trajectories are proven in \cite[Section VII]{3Diekmann}.
We here merely will use the result on invariance of the $\ome$-limit set - for an alternative proof of Corollary \ref{corol1} below. 

\begin{lemma}\label{lem20}
Let $\s:I\-> X$ be a trajectory, then $\s$ is continuous. If $\sup I=\infty$, then
\[
\ome(\s)=\bigcap_{t\ge 0} \overline{\s(I\cap[t,\infty))}.
\]
If additionally $\overline{\s(I)}$ is compact, then $\ome(\s)$ is nonempty, compact and connected, $dist(\s(t),\ome(\s))\->0$ as $t\ri>\infty$ and for $x\in\ome(\s)$ one has $\Si(t,x)\in\ome(\s)$ for all $t\ge 0$. 
\end{lemma}
We now conclude continuity of the semiflow from continuous 
dependence on initial values and the semigroup property from uniqueness. 
In the following we assume that $t_0=0$. 
\begin{theorem}\label{theo8}
Suppose that $f:\cD\->\R^n$ is almost locally Lipschitz and satisfies (F),
that $f(U)$ is bounded whenever $U\subset\cD$ is bounded and that $\{x_t^\phi:\;t\in[0,\a)\}$ is bounded whenever $\phi\in V_\cD$ and whenever $x^\phi$ is defined on $[0,\a)$. Then for any $\phi\in V_\cD$ there exists a unique global solution and for all $t\ge0$ one has $x_t^\phi\in V_\cD$. Hence,  we can define a map 
\[
S:[0,\infty)\times V_\cD\->V_\cD;\;S(t,\phi):=x_t^\phi.
\]
This map defines a continuous semiflow on $V_\cD$ with respect to the $\sup$-norm. 
\end{theorem}
\Proof Existence of a unique global solution for all $\phi\in V_\cD$ follows from Theorem \ref{theo6}. 
Let $\phi\in V_\cD$ and $t>0$. By definition of a solution we have $x_t^\phi\in\cD$. Let $r,s\in[-h,
0]$. Then
\[
|x^\phi_t(r)-x^\phi_t(s)|=|x^\phi(t+s)-x^\phi(t+r)|. 
\]
First, $x^\phi$ is Lipschitz on $[-h,0]$, since $\phi$ is Lipschitz. Next, $x^\phi$ is as a solution 
differentiable on $[0,t]$ and satisfies (\ref{eq7}--\ref{eq1}). Hence, $(x^\phi)'$ is continuous. Thus $x^\phi$ is Lipschitz on $[0,t]$ by the mean value theorem. Hence $x^\phi$ is Lipschitz on $[-h,t]$ and thus $x_t^\phi\in V_\cD$. Next, it is clear that $S(0,\phi)=\phi$ for all
$\phi\in V_\cD$. 

To see that the semigroup property holds, fix $\phi$ and define for some $t>0$ and $\t>0$
\begin{eqnarray}
y(s)&:=&\begin{cases}
\phi(s),&s\in[-h,0]
\\
S(s,\phi)(0),& s\in[0,t]
\\
S(s-t,S(t,\phi))(0),&s\in[t,t+\t]
\end{cases}
\nn\\
z(s)&:=&\begin{cases}
\phi(s),&s\in[-h,0]
\\
S(s,\phi)(0),& s\in[0,t+\t].
\end{cases}
\nn
\end{eqnarray}
We have $y=z$ on $[-h,t]$, hence in particular on $[t-h,t]$, thus $y_t=z_t$. 
Now suppose that 
$s\in(t,t+\t]$. Let $\th\in[-h,0]$. If $s-t+\th\ge0$, then 
\begin{eqnarray}
x_{s-t}^{S(t,\phi)}(\th)=x_{s-t+\th}^{S(t,\phi)}(0)=S(s+\th-t,S(t,\phi))(0)=y(s+\th)=y_s(\th).
\nn
\end{eqnarray}
If $s-t+\th<0$ then
\begin{eqnarray}
x_{s-t}^{S(t,\phi)}(\th)&=&S(t,\phi)(s-t+\th)=x_t^{\phi}(s-t+\th)=x^\phi(s+\th)
\nn\\
&=&\begin{cases}
\phi(s+\th),&s+\th\le 0
\\
S(s+\th,\phi)(0),&s+\th>0
\end{cases}
=y(s+\th)=y_s(\th). 
\nn
\end{eqnarray}
Thus $x_{s-t}^{S(t,\phi)}=y_s$. Hence
\begin{eqnarray}
y'(s)=\frac{d}{ds}x_{s-t}^{S(t,\phi)}(0)=(x^{S(t,\phi)})'(s-t)=f(x_{s-t}^{S(t,\phi)})=f(y_s). 
\nn
\end{eqnarray}
Hence  with $t$ and $t_0$ replaced by $s$ and $t$ respectively, $y$ is a solution 
of (\ref{eq7}) through $z_t$ on $[-h,t+\t]$. One similarly shows that so is $z$. By uniqueness we have $y=z$ on $[-h,t+\t]$. If we fill $s=t+\t$ and use the
definitions of $y$ and $z$, we see that this implies the semigroup property.

To see continuity of $S$ note that 
\begin{eqnarray}
&&|S(t,\phi)(\th)-S(\ote,\ophi)(\th)|
\nn\\
&\le&|S(t,\phi)(\th)-S(t,\ophi)(\th)| +|S(t,\ophi)(\th)-S(\ote,\ophi)(\th)|
\nn\\
&=&|x^\phi(t+\th)-x^\ophi(t+\th)| +|x^\ophi(t+\th)-x^\ophi(\ote+\th)|.
\nn
\end{eqnarray}
The first term can be estimated using our result on continuous 
dependence (Theorem \ref{theo7}), the second using continuity of solutions in time. 
\qed

Continuous dependence and the semigroup property can be combined to prove the following
result:
\begin{corol}\label{corol1}
Suppose that $f$ satisfies the assumptions of Theorem \ref{theo8},
$\phi\in V_\cD$, $x^\phi(t)\->x^*\in\R$ as $t\ri>\infty$. Then $x^*$ is an equilibrium solution. 
\end{corol}
\Proof Let $(t_k)\in[0,\infty)^\N$, $t_k\ri>\infty$, and fix $t>0$. Define a sequence via $\phi^k:=S(t_k,\phi)\in\cD$ and
denote by $\phi^*$ the constant function with value $x^*$ on $[-h,0]$. Then 
$\phi^k=x^\phi_{t_k}\->\phi^*$ (uniformly) by our assumption. Similarly $S(t+t_k,\phi)\->\phi^*$. 
But also $S(t+t_k,\phi)=S(t,S(t_k,\phi))=S(t,\phi^k)\->S(t,\phi^*)$ by Theorem \ref{theo8}. Hence $S(t,\phi^*)=\phi^*$. One can conclude that $x^*$ is an 
equilibrium solution. 
\qed

The result can also be concluded from Lemma \ref{lem20}:
\bigskip

\noindent
{\bf Proof of Corollary \ref{corol1} via Lemma \ref{lem20}.} 
Define $I:=[0,\infty)$, choose any $\phi\in V_\cD$, and define $\s(t):=S(t,\phi)$. Then
$\s$ is a trajectory. We show that $\overline {\s(I)}$
is compact, i.e., that $\sigma(I)=S([0,\infty),\phi)$ is relative compact. Let $(t_n)\in
[0,\infty)^\N$. Case 1: $t_n\in[0,T]$ for all $n\in\N$ and some $T>0$. Hence, there exists
$(t_{n_j})\subset (t_n)$, $\ote\in[0,T]$ such that $t_{n_j}\->\ote$ as $j\ri>\infty$. Then 
$(S(t_{n_j},\phi))\subset (S(t_n,\phi))$ and $S(t_{n_j},\phi)\->S(\ote,\phi)$ by continuity.
Case 2: $(t_n)$ is unbounded. Then there exists some $(t_{n_j})\subset(t_n)$ such that
$t_{n_j}\->\infty$ as $j\ri>\infty$. Thus $S(t_{n_j},\phi)\->\phi^*$ where $\phi^*\in V_\cD$
is defined as $\phi^*(t)=x^*$ for all $t\in[-h,0]$. Hence, in any case, 
$(S(t_n,\phi))$ has a Cauchy subsequence, thus $\overline {\s(I)}$ is compact.  Now
note that for the $\ome$-limit set of the trajectory one has $\ome(\s)=\{\phi^*\}$. Then by Lemma \ref{lem20} one has $S(t,\phi^*)\in\ome(\s)=\{\phi^*\}$, 
i.e., $S(t,\phi^*)=\phi^*$ for all $t\ge 0$. 
\qed

%
%
%
\section{Invariant compact sets}\label{s3}
\subsection{Assumptions, main results and discussion}
In the setting of Section \ref{s2} we now set $n=2$ and $B=0$, such that $\cD=C([-h,0],\R_+^2)\subset \cC=C([-h,0],\R^2)$, where $\R_+=[0,\infty)$, and consider a functional
$j:\cD\->\R_+$,  a function $q:\R_+\->\R$ and a DDE of the form 
\begin{eqnarray}
w'=q(v)w,\;\;v'(t)=-\mu v(t)+j(w_t,v_t),\;\;t>0,\;\;(w_0,v_0)=(\vi,\psi)\in\cD,
\nn\\
\label{eq18}
\end{eqnarray}
where $\mu>0$ is a parameter. 
Define $\oq:=\sup q$ and suppose throughout the section that $\oq<\infty$, $q$ is locally Lipschitz, $j$ is almost
locally Lipschitz and that for some $k_j>0$ at least one of the two, 
\begin{eqnarray}
j(\vi,\psi)\le k_j\|\vi\|,\;\;{\rm or}
\label{eq19}\\
j(\vi,\psi)\le k_j\vi(-\t(\psi)),
\label{eq20}
\end{eqnarray}
where $\t:C([-h,0],\R_+)\->[\uta,h)$ for some $\uta\in(0,h)$, holds. 

Obviously (\ref{eq19}) is a weaker requirement. As we will see, however, 
(\ref{eq20}) may lead to better results while still applicable to our model. 
Our proofs in the context of invariant sets of bounded functions rely on an exponential
estimate for the $w$-component that uses the linearity of the $w$-equation. Exponential estimates can be derived for general DDE, see e.g. \cite[Corollary 6.1.1]{Hale}, so our approach possibly works for systems more general than (\ref{eq18}) too. In the context of our application, however, we found (\ref{eq18}) a good compromise between the wishes to be general and to provide sharp estimates for our model.

Now note that, supposing a solution through $(\vi,\psi)\in\cD$ exists,  one has that 
\begin{eqnarray}
w(t)&=&\begin{cases}
\vi(t),&t\in[-h,0]
\\
\vi(0)e^{\int_0^tq(v(s))ds},&t>0,
\end{cases}
\label{eq21}
\end{eqnarray}
hence
\begin{eqnarray}
w(t)&\le&\|\vi\|q_e(t), \forall\;t\ge-h, 
\;\;{\rm where}\;\;
q_e(t):=\begin{cases}
1,&t\in[-h,0]
\\
e^{\oq t},&t>0. 
\end{cases}
\label{eq22}
\end{eqnarray}
Note that $q_e$ is continuous, nondecreasing, increasing on $[0,\infty)$ and differentiable
on $[-h,0)\cup(0,\infty)$.

An important case is that $q$ is decreasing and has one positive zero, see also Section 
\ref{ss5} and \cite{Getto1}. Hence, positivity of $q$ is not out-ruled, and 
thus, looking at (\ref{eq21}), we 
cannot expect that a next-state operator $\phi\mapsto S(t,\phi)$  maps
a set of the form $C([-h,0],[0,A]\times[0,B])$, $A,B\in(0,\infty)$ into itself
(to avoid subindices, we here, other than in the previous section, let both $A$ and $B$ denote upper bounds). For similar reasons (see the proof of Theorem \ref{theo10}
(b) below) we cannot expect this for a set of $R$-Lipschitz functions either. On the other hand, filling (\ref{eq21})
into the second equation of (\ref{eq18}) yields a closed system in $v$ (depending on both initial histories). 
Motivated by this, next to an initial result for the $w$-component,  we will establish an invariant set for the $v$-component. We refer to the discussion section for possible extensions of this research. 

Define for any $B>0$ and $R>0$ the set 
\begin{eqnarray}
C_{B,R}:=\{\chi\in C([-h,0],[0,B]),\;lip\;\chi\le R\}.
\label{eq25}
\end{eqnarray}
Note that $C_{B,R}$ is convex and, by the Arzela-Ascoli theorem, compact.  
Next, we formulate the main
results of this section and give proofs in the next subsection. 
With the cases (\ref{eq19}--\ref{eq20}), respectively, we associate functions $f_l,f_\tau:\R_+\->\R_+;$
\begin{eqnarray}
f_l(t)&:=&\frac{k_j}{\mu+\oq}(e^{\oq t}-e^{-\mu t}),
\nn\\
f_\t(t)&:=&
\begin{cases}
\frac{k_j}{\mu}(1-e^{-\mu t}),&{\rm if}\;t\le\uta
\\
k_j\frac{\oq(e^{-\mu(t-\uta)}-e^{-\mu t})+\mu(e^{\oq (t-\uta)}-e^{-\mu t})}{\mu(\mu+q)},&{\rm if}\;t>\uta,
\end{cases}
\nn
\end{eqnarray}
(where l stands for linear in reference to (\ref{eq19})). When writing about
these functions we will assume that the respective case holds, sometimes only implicitly.  
\begin{theorem}\label{theo10}
Under the assumptions of this subsection, the following holds for any $(\vi,\psi)\in V_\cD$ . 
\begin{itemize}
\item[(a)]  The system (\ref{eq18}) has a unique solution $x=(w,v)$ through $(\vi,\psi)$ on 
$[0,\infty)$. The solutions define a continuous semiflow in the sense
of Theorem \ref{theo8}.
\item[(b)] Choose $A$, $R$ and $T$ such that $\oq A e^{\oq T}\le R$. Then, if $\|\vi\|\le A$ and $lip\;\vi\le R$ one has $lip\;w_t\le R$ for all $t\in[0,T]$. 
\item[(c)] Both, $f_l$ and $f_\tau$ are
zero in zero, tend to $\infty$ at $\infty$, are increasing and continuous, $f_l$ is differentiable, and $f_\tau$ is differentiable on $[0,\uta)\cup(\uta,\infty)$. The functions
\[
t\longmapsto\frac{f_l(t)}{1-e^{-\mu t}}\;{\rm and}\;
t\longmapsto\frac{f_\tau(t)}{1-e^{-\mu t}},
\]
respectively, increase from $k_j/\mu$ to infinity on $\R_+$, and equal $k_j/\mu$ on $[0,\uta]$ and increase to infinity on $[\uta,\infty)$. Finally, $f_l(t)>f_\t(t)$ for all $t>0$. 
\item[(d)]
Assume that (\ref{eq19}) holds and choose $A$, $B$, $R$ and $T$ such that 
$\frac{Af_l(T)}{1-e^{-\mu T}}\le B$ and $R\ge \max\{\mu B, k_jAe^{\oq T}\}$. Then, if $\|\vi\|\le A$ and $\psi\in C_{B,R}$ one has $v_t\in C_{B,R}$ for all $t\in[0,T]$. 
\end{itemize}
If (\ref{eq20}) holds, then the following hold. 
\begin{itemize}

\item[(e)]
Choose $A$, $B$, $R$ and $T$ such that $\frac{Af_\tau(T)}{1-e^{-\mu T}}\le B$ and 
\[
R\ge \max\{\mu B, k_jAq_e(T-\uta)\}.
\] 
Then, if $\|\vi\|\le A$ and $\psi\in C_{B,R}$, one has $v_t\in C_{B,R}$ for all $t\in[0,T]$. 

\item[(f)]
Choose $A$, $B$ and $R$ such that 
$Ak_j< B\mu\le R$ and $\de$ such that
$Ak_je^{\oq\de}=\mu B$. 
Then, if $\|\vi\|\le A$ and $\psi\in C_{B,R}$, one has $v_t\in C_{B,R}$ for all $t\in[0,\uta+\de]$. 

\end{itemize}
\end{theorem}
For further discussion of the theorem we state some technical results. 
\begin{lemma}\label{lem31}
One has $\frac{k_je^{\oq t}}{\mu}>\frac{f_l(t)}{1-e^{-\mu t}}$ for all $t>0$
and $k_je^{\oq(t-\uta)}>\mu\frac{f_\t(t)}{1-e^{-\mu t}}$ for $t>\uta$.
\end{lemma}
Now, note that (f) is a simple corollary of (e). To prove this, define $T=\uta+\de$ in (e), and apply the second estimate of the lemma with $t=T$.  We omit further details. 

By the previous lemma, in Theorem \ref{theo10} (d) and (e) it would be sufficient to assume
that 
\begin{eqnarray}
R\ge\mu B\ge
\begin{cases}
k_jAe^{\oq T},& {\rm respectively},
\\
k_jAe^{\oq(T-\uta)},&
\nn
\end{cases}
\nn
\end{eqnarray}
which is stronger but easier to check than the present assumptions. 

Note that (e) allows to establish for the solution a lower bound and a lower Lipschitz constant
than (d): Fix $A$ and $T$. Then the lowest bound we can achieve through (d) is $B_d:=\frac{Af_l(T)}{1-e^{-\mu T}}$, whereas through (e) we can achieve the bound $B_e:=\frac{Af_\t(T)}{1-e^{-\mu T}}<B_d$. The lowest Lipschitz constant we can achieve through (d) is 
$R_d:=\max\{\mu B_d,k_jAe^{\oq T}\}>\max\{\mu B_e,k_jAq_e(T-\uta)\}=:R_e$, where $R_e$ is a
(the lowest) Lipschitz constant we can achieve through (e). 

We get invariance for a longer time through (e) than through (d): Fix $A$, $B$ and $R$ such that $\frac{Ak_j}{\mu}<B$ and $R\ge\mu B$. Then the largest time spans which (d) and (e) yield are respectively $t_d:=\min\{t_{d1},t_{d2}\}$ and $t_e=\min\{t_{e1},t_{e2}\}$, where the involved quantities are defined via $\frac{Af_l(t_{d1})}{1-e^{-\mu t_{d1}}}=B$, $\frac{Af_\t(t_{e1})}{1-e^{-\mu t_{e1}}}=B$, $R=\max\{\mu B,Ak_je^{\oq t_{d2}}\}$
and $R=\max\{\mu B,Ak_jq_e(t_{e2}-\uta)\}$. One has $t_{dj}<t_{ej}$, $j=1,2$, 
hence $t_d<t_e$. 

Theorem \ref{theo10} (f) shows that, if (\ref{eq20}) holds,
there is a lower bound ($\uta$) for the time for which invariance holds, which is uniform for  all $A$, $B$ satisfying $\frac{Ak_j}{\mu}<B$. If merely (\ref{eq19})
holds we cannot get such a lower bound through (d). 
%
%
%
\subsection{Proofs}
We start with some general facts regarding the computation with almost locally Lipschitz
functions. In the following lemma we let $f$  and $g$ denote arbitrary functions and $\cD\subset \cC$ an arbitrary domain. 
\begin{lemma}\label{lem16} 
(a) Suppose that $f,g:\cD\subset \cC\->\R$ are almost locally Lipschitz. Then so are $fg$, $(f,g)$ and $f+g$. 
\\
(b) Let  $f:\cD\subset\cC\->\R$ be almost locally Lipschitz and $g:f(\cD)\subset\R\->\R$
locally Lipschitz, then $g\circ f:\cD\->\R$ is almost locally Lipschitz. 
\end{lemma}
\Proof (a) Clearly $fg$ is continuous. Now let $\phi_0\in\cD$, $R>0$. Choose $\de_{f}$, $k_{f}$, $\de_{g}$, $k_{g}$
in notation similar as in Definition \ref{def1} and according to the definition. Define $k:=\max\{
k_{f},k_{g}\}$. By continuity of $f$ and $g$ we can 
choose $M$ and $\de_1$ such that $f$ and $g$ are bounded by $M$ on $\oB_{\de_1}(\phi_0)$. Define 
$\de:=\min\{\de_{f},\de_{g},\de_1\}$. Then $fg$ is $k$-Lipschitz on $V(\phi_0;R,\de)$, hence almost locally Lipschitz. The remainder of the proof of (a) is obvious. 
\\
(b) First, clearly $g\circ f$ is continuous. Next, let $\phi_0\in\cD$, $R>0$, choose $\e$, $k_1$ such that $g$ is $k_1$-Lipschitz 
on $B_\e(f(\phi_0))$. Choose, $\de$, $k_2$ such that $f$ is $k_2$-Lipschitz on 
$V(\phi_0;\de,R)$ and $f(B_\de(\phi_0))\subset B_\e(f(\phi_0))$. Let $\vi, \psi\in V(\phi_0;\de,R)$.
Then the following estimate implies the statement:
\begin{eqnarray}
|g(f(\vi))-g(f(\psi))|\le k_1|f(\vi)-f(\psi)|\le k_1k_2\|\vi-\psi\|. 
\nn
\end{eqnarray}
\qed

In view of applying the general theory we next would like to show that a functional $f$ 
associated with (\ref{eq18}) is almost locally Lipschitz. 
Due to the previous result it is sufficient to show that so are the components of $f$  and we 
start with the first component. 
\begin{lemma}\label{lem21} The functional $f_1:\cD\->\R;\; f_1(\vi,\psi):=q(\psi(0))\vi(0)$ is locally Lipschitz, in particular 
almost locally Lipschitz. 
\end{lemma}
\Proof First note that the projection map and the evaluation map 
\begin{eqnarray}
&&C([-h,0],\R_+^2)\->C([-h,0],\R_+);\;(\vi,\psi)\mapsto\vi,
\;{\rm and}
\nn\\
&&C([-h,0],\R_+)\->\R;\;\vi\mapsto\vi(0)
\nn
\end{eqnarray}
and analogous maps for the $\psi$-component are locally Lipschitz. Hence, by the
preservation of local Lipschitzianity under composition and the Lipschitz property of $q$
it follows that $(\vi,\psi)\mapsto q(\psi(0))$ is locally Lipschitz. Moreover $(\vi,\psi)\mapsto\vi(0)$ is locally Lipschitz. Thus by the product rule for locally Lipschitz functions so is $f_1$. 
\qed

\noindent
{\bf Proof of Theorem \ref{theo10} (a).} By Lemmas \ref{lem16} and \ref{lem21} it follows that 
\begin{eqnarray}
f(\vi,\psi)=(q(\psi(0))\vi(0),-\mu\psi(0)+j(\vi,\psi))^T
\nn
\end{eqnarray}
is almost locally Lipschitz. Property $(F)$ is guaranteed by non-negativity of $j$. The 
boundedness property of $f$ required in Theorem \ref{theo8} is guaranteed by continuity of
$q$
and (\ref{eq19}--\ref{eq20}). The required boundedness property of the trajectory can be
guaranteed by (\ref{eq21}--\ref{eq22}) if one integrates the $v$-equation in (\ref{eq18}) using
$\oq<\infty$ and the
variations of constants formula. Application of Theorem \ref{theo8} completes the proof. \qed

\bigskip

\noindent
{\bf Proof of Theorem \ref{theo10} (b).}  Let $t\in[0,T]$. It is equivalent to show that $lip\;w|_{[t-h,t]}\le R$. Since $lip\;\vi\le R$ it follows that 
$w$ is $R$-Lipschitz on $[t-h,t]\cap[-h,0]$. On $[t-h,t]\cap[0,\infty)$ the function $w$ is 
differentiable with
\begin{eqnarray}
|w'(t)|\le|q(v(t))|\|\vi\|q_e(t)\le\oq A e^{\oq t}\le\oq A e^{\oq T}\le R.
\nn
\end{eqnarray}
\qed

In subsequent proofs we will sometimes omit bars in $\oq$ and $\uta$ for the sake of the presentation. 
\begin{lemma}\label{lem22}
One has for any $t>0$
\begin{eqnarray}
v(t)&\le&
\begin{cases}
e^{-\mu t} \psi(0)+\|\vi\|f_l(t),&{\rm if} (\ref{eq19})\;{\rm holds}, 
\\
e^{-\mu t} \psi(0)+\|\vi\|f_\tau(t),&{\rm if} (\ref{eq20})\;{\rm holds}. 
\end{cases}
\label{eq23}
\end{eqnarray}
Moreover $f_l(t)>f_\t(t)$ for all $t>0$. 
\end{lemma}
\Proof By the variation of constants formula 
\begin{eqnarray}
v(t)=e^{-\mu t}\psi(0)+e^{-\mu t}\int_0^te^{\mu s}j(w_s,v_s)ds. 
\nn
\end{eqnarray}
If (\ref{eq19}) holds, 
\begin{eqnarray}
e^{-\mu t}\int_0^te^{\mu s}j(w_s,v_s)ds\le k_je^{-\mu t}\int_0^te^{\mu s}\|w_s\|ds\le
\|\vi\|k_je^{-\mu t}\int_0^te^{(\mu+q)s}ds, 
\nn
\end{eqnarray}
which yields the first statement. 
If (\ref{eq20}) holds, then
\begin{eqnarray}
&&e^{-\mu t}\int_0^te^{\mu s}j(w_s,v_s)ds\le e^{-\mu t}k_j\int_0^te^{\mu s}w(s-\t(v_s))ds
\nn\\
&\le&e^{-\mu t}k_j\|\vi\|\int_0^te^{\mu s}q_e(s-\t(v_s))ds
\le e^{-\mu t}k_j\|\vi\|\int_0^te^{\mu s}q_e(s-\t)ds.
\nn
\end{eqnarray}
If $t\le\t$ the statement is obvious. If $t>\t$, then 
\begin{eqnarray}
&&e^{-\mu t}\int_0^te^{\mu s}j(w_s,v_s)ds\le k_j\|\vi\|
[e^{-\mu t}\int_0^\t e^{\mu s}ds+e^{-\mu t}\int_\t^te^{q(s-\t)+\mu s}ds],
\nn
\end{eqnarray}
which also yields the first statement. Now note that by the above estimates
\begin{eqnarray}
f_l(t)=k_je^{-\mu t}\int_0^te^{(\mu+q)s}ds, \;\;
f_\t(t)=k_je^{-\mu t}\int_0^te^{\mu s}q_e(s-\t)ds. 
\nn
\end{eqnarray}
Hence $f_l(t)>f_\t(t)$ for all $t>0$ if $e^{qs}>q_e(s-\t)$ for all $s>0$, which is the case. 
\qed

\noindent
{\bf Proof of Theorem \ref{theo10} (c).} First note that by the previous lemma
$f_l(t)>f_\t(t)$ for all $t>0$. 
Next, if (\ref{eq19}) holds,
\begin{eqnarray}
&&{\rm sgn}\frac{d}{dt}\frac{f_l(t)}{1-e^{-\mu t}}
={\rm sgn}[(qe^{qt}+\mu e^{-\mu t})(1-e^{-\mu t})-(e^{qt}-e^{-\mu t})\mu e^{-\mu t}]
\nn\\
&=&{\rm sgn}[qe^{qt}-(q+\mu)e^{(q-\mu)t}+\mu e^{-\mu t}]={\rm sgn}\;g(t)
\nn
\end{eqnarray}
for obviously defined $g$. Then $g(0)=0$ and 
\begin{eqnarray}
g'(t)&=&q^2e^{qt}+(\mu^2-q^2)e^{(q-\mu)t}-\mu^2e^{-\mu t}
\nn\\
&=&q^2e^{qt}(1-e^{-\mu t})+\mu^2e^{-\mu t}(e^{qt}-1)>0.
\nn
\end{eqnarray}
Thus $g(t)>0$ for all $t>0$ and hence $t\mapsto f_l(t)/(1-e^{-\mu t})$ is increasing. 

If (\ref{eq20}) holds, to see that $t\mapsto f_\tau(t)/(1-e^{-\mu t})$ is increasing, it is sufficient to show that
\begin{eqnarray}
g(t):=\frac{q(e^{-\mu(t-\t)}-e^{-\mu t})+\mu(e^{q (t-\t)}-e^{-\mu t})}{1-e^{-\mu t}}
\nn
\end{eqnarray}
is increasing for $t>\t$. One has
\begin{eqnarray}
&&{\rm sgn}\;g'(t)={\rm sgn}\;\{[q(\mu e^{-\mu t}-\mu e^{-\mu(t-\t)})+\mu(qe^{q(t-\t)}+\mu e^{-\mu t})](1-e^{-\mu t})
\nn\\
&&-\mu e^{-\mu t}[q(e^{-\mu(t-\t)}-e^{-\mu t})+\mu(e^{q(t-\t)}-e^{-\mu t})]\}
\nn\\
&=&{\rm sgn}\;\{[q(e^{-\mu t}- e^{-\mu(t-\t)})+qe^{q(t-\t)}+\mu e^{-\mu t}](1-e^{-\mu t})
\nn\\
&&
-e^{-\mu t}[q(e^{-\mu(t-\t)}-e^{-\mu t})+\mu(e^{q(t-\t)}-e^{-\mu t})]\}
\nn\\
&=&{\rm sgn}\;\{q(e^{-\mu t}- e^{-\mu(t-\t)})+qe^{q(t-\t)}+\mu e^{-\mu t}
+q(e^{-\mu(2t-\t)}-e^{-2\mu t})
\nn\\
&&-qe^{(q-\mu)t-q\t}-\mu e^{-2\mu t}+q(e^{-2\mu t}-e^{-\mu(2t-\t)})
+\mu(e^{-2\mu t}-e^{(q-\mu)t-q\t})\}
\nn\\
&=&{\rm sgn}\;\{q(e^{-\mu t}- e^{-\mu(t-\t)})+qe^{q(t-\t)}+\mu e^{-\mu t}
-(q+\mu)e^{(q-\mu)t-q\t}\}
\nn\\
&=&sgn\;h(q)
\nn
\end{eqnarray}
for obviously defined $h$. Then $h(0)=0$. Next, 
\begin{eqnarray}
h'(q)&=&e^{-\mu t}- e^{-\mu(t-\t)}+e^{q(t-\t)}+q(t-\t)e^{q(t-\t)}
\nn\\
&&-[e^{(q-\mu)t-q\t}+(q+\mu)(t-\t)e^{(q-\mu)t-q\t}],
\nn\\
h'(0)&=&1-e^{-\mu(t-\t)}-\mu(t-\t)e^{-\mu t}=:j(t)
\nn
\end{eqnarray}
in obvious notation. Then $j'(t)=\mu e^{-\mu t}[e^{\mu\t}-1+\mu(t-\t)]>0$, hence $j(t)> j(\t)=0$
and thus $h'(0)> 0$. Next, 
\begin{eqnarray}
h''(q)&=&(t-\t) e^{q(t-\t)}\{2+q (t-\t)-[2+(q+\mu)(t-\t)]e^{-\mu t}\}
\nn\\
&=&(t-\t)e^{q(t-\t)}k(q)
\nn
\end{eqnarray}
for obviously defined $k$. Then,  applying $e^x\ge 1+x$ to $x=\mu(t-\t)$,
\begin{eqnarray}
k(0)&=&2-[2+\mu (t-\t)]e^{-\mu t}\ge 1-e^{-\mu t}+1-e^{-\mu\t}> 0, 
\nn\\
k'(q)&=&t-\t-e^{-\mu t}(t-\t)>0.
\nn
\end{eqnarray}
Hence, $k$ is positive for $q>0$, thus so is $h''$, hence so is $h'$, thus so is $h$, 
hence so is $sgn\;g'$. We have shown that $t\mapsto f_\t(t)(1-e^{-\mu t})$ is 
increasing. Monotonicity of $f_l$ follows from monotonicity of 
$f_l(t)/(1-e^{-\mu t})$ and the same conclusion holds for $f_\tau$. 
Using that $(1-e^{-\mu t})^{-1}$ is bounded at infinity the remaining statements are
easy to see. 
\qed

\begin{lemma}\label{lem28}
Assume that (\ref{eq19}) holds and that $A$, $B$ and $T$ are such that 
$\frac{Af_l(T)}{1-e^{-\mu T}}\le B$. Then $\|\vi\|\le A$ and $\|\psi\|\le B$ imply
that $v(t)\le B$ for all $t\in[-h,T]$. 
\end{lemma}
\Proof By (\ref{eq23}) one has 
$
v(t)\le Be^{-\mu t}+Af_l(t)
$
for $t\in(0,T]$.
Hence $v(t)\le B$ if $Af_l(t)/(1-e^{-\mu t})\le B$ and the latter follows by assumption and Theorem \ref{theo10} (c). 
\qed

An elaboration of the maximum in the following lemma will be carried out further down. 
\begin{lemma}\label{lem27}
Let $\|\vi\|\le A$ and $\|\psi\|\le B$. Let $T>0$ and choose
\begin{eqnarray}
R&\ge&
\begin{cases}
\max_{t\in[T-h,T]\cap[0,\infty)}\max\{k_jq_e(t)A, \mu(e^{-\mu t}B+Af_l(t))\},&
\\
{\rm if}\;(\ref{eq19})\; {\rm holds},&
\\
\max_{t\in[T-h,T]\cap[0,\infty)}\max\{k_jq_e(t-\uta)A, \mu(e^{-\mu t}B+Af_\tau(t))\},
\\
{\rm if}\;(\ref{eq20})\; {\rm holds}.
\end{cases}
\nn
\end{eqnarray}
Then, if $lip\;\psi\le R$, also $lip\;v_T\le R$. 
\end{lemma}
\Proof  We should show that $lip\;v|_{[T-h,T]}\le R$. First, 
\begin{eqnarray}
lip\;v|_{[T-h,T]\cap[-h,0]}=lip\;\psi|_{[T-h,T]\cap[-h,0]}
\le R. 
\nn
\end{eqnarray}
Next, if (\ref{eq19}) holds, we get
$
v'(t)\le j(w_t,v_t)\le k_j\|w_t\|\le k_j q_e(t)\|\vi\|. 
$
If (\ref{eq20}) holds, then
$
v'(t)\le  k_jw(t-\t(v_t))\le  k_j q_e(t-\t(v_t))\|\vi\|\le k_j q_e(t-\uta)\|\vi\|. 
$
Moreover
\begin{eqnarray}
v'(t)\ge-\mu v(t)\ge \begin{cases}
-\mu(e^{-\mu t}|\psi(0)|+\|\vi\|f_l(t)),&{\rm if}\;(\ref{eq19})\;{\rm holds}
\\
-\mu(e^{-\mu t}|\psi(0)|+\|\vi\|f_\tau(t)),&{\rm if}\;(\ref{eq20})\;{\rm 
holds}.
\end{cases}
\nn
\end{eqnarray}
Hence for $t>0$ one has
\begin{eqnarray}
|v'(t)|&\le&
\max\{k_jq_e(t)\|\vi\|,\mu(|\psi(0)|e^{-\mu t}+\|\vi\|f_l(t))\},
\nn\\
&&{\rm if}\;(\ref{eq19})\;{\rm holds}
\nn\\
|v'(t)|&\le&\max\{k_jq_e(t-\uta)\|\vi\|,\mu(|\psi(0)|e^{-\mu t}+\|\vi\|f_\tau(t))\},
\nn\\
&&{\rm if}\;(\ref{eq20})\;{\rm 
holds}.
\nn
\end{eqnarray}
Hence $
lip\;v|_{[T-h,T]\cap[0,\infty)}\le\max_{t\in[T-h,T]\cap[0,\infty)}|v'(t)|\le R.
$
\qed

\begin{lemma}\label{lem26}
Assume that (\ref{eq19}) holds, choose $A$ and $B$ such that $Ak_j/\mu<B$ and define $t_1$ via $A\frac{f_l(t_1)}{1-e^{-\mu t_1}}=B$, then if $T\le t_1$
one has
\begin{eqnarray}
\max_{t\in[0,T]}\mu[ e^{-\mu t}B+Af_A(t)]=\mu B. 
\nn
\end{eqnarray}
\end{lemma}
\Proof  Define $g:[0,t_1]\->\R_+;\;g(t):=\mu[ e^{-\mu t}B+Af_l(t)]$. Then,
\begin{eqnarray}
g'(t)&=&\mu[\frac{Ak_j}{\mu+q}(qe^{qt}+\mu e^{-\mu t})-\mu Be^{-\mu t}], \;
g'(0)=\mu(Ak_j-\mu B)<0, 
\nn\\
g'(t_1)&=&B\mu[\frac{(1-e^{-\mu t_1})(qe^{qt_1}+\mu e^{-\mu t_1})}{e^{q t_1}-e^{-\mu t_1}}
-\mu e^{-\mu t_1}]
\nn\\
&=&
\frac{B\mu}{e^{q t_1}-e^{-\mu t_1}}
[qe^{qt_1}+\mu e^{-\mu t_1}-(q+\mu )e^{-(\mu-q) t_1}]
=\frac{B\mu}{e^{q t_1}-e^{-\mu t_1}}h(t_1)
\nn
\end{eqnarray}
for obviously defined $h$. We have seen in the proof of Theorem \ref{theo10} (c) that $h(t_1)>0$. 
Thus $g'(t_1)>0$. Next

\begin{eqnarray}
g''(t)&=&\mu[\frac{Ak_j}{\mu+q}(q^2e^{qt}-\mu^2e^{-\mu t})+\mu^2Be^{-\mu t}]
\nn\\
&=&\mu[\frac{Ak_j}{\mu+q}(q^2e^{qt}-\mu^2e^{-\mu t})+\mu^2e^{-\mu t}\frac{Af_l(t_1)}{1-e^{-\mu t_1}}]
\nn\\
&>&\mu[\frac{Ak_j}{\mu+q}(q^2e^{qt}-\mu^2e^{-\mu t})+\mu^2e^{-\mu t}\frac{Af_l(t)}{1-e^{-\mu t}}]
\nn\\
&=&\frac{\mu Ak_j}{\mu+q}[q^2e^{qt}-\mu^2e^{-\mu t}+\mu^2e^{-\mu t}\frac{e^{qt}-e^{-\mu t}}{1-e^{-\mu t}}]
\nn\\
&=&\frac{\mu Ak_j}{(\mu+q)(1-e^{-\mu t})}[q^2e^{qt}(1-e^{-\mu t})+\mu^2e^{-\mu t}
(e^{qt}-1)]>0. 
\nn
\end{eqnarray}
Hence, $g'$ increases monotonously from a negative value to a positive value. Thus $g$
decreases monotonously to a minimum, then increases monotonously, hence assumes a
maximum either in zero, or in $t_1$. Since $g(0)=g(t_1)=\mu B$ the statement follows. 
\qed

\noindent
{\bf Proof of Theorem \ref{theo10} (d).}
 First note that $T<t_1$ for $t_1$ as in Lemma \ref{lem26}. Hence by this lemma and
Lemma \ref{lem27} one has $lip\;v_t\le R$ for all $t\in[0,T]$. The boundedness property follows by
Lemma \ref{lem28}. 
\qed

\bigskip

\noindent
{\bf Proof of Theorem \ref{theo10} (e).} The stated boundedness is implied by the monotonicity shown in (c) and Lemma \ref{lem22}. 
Moreover, since $\mu B\ge\mu( Be^{-\mu t}+Af_\tau(t))$ for $t\in[0,T]$, one has
\begin{eqnarray}
R&\ge&\max\{Ak_jq_e(T-\t),\mu B\}
\nn\\
&\ge&\max_{t\in[0,T]}\max\{Ak_jq_e(t-\t),\mu( Be^{-\mu t}+Af_\tau(t))\}.
\nn
\end{eqnarray}
Hence the Lipschitz-property follows by Lemma \ref{lem27}.
\qed

\bigskip

\noindent
{\bf Proof of Lemma \ref{lem31}.} For $t>s\ge0$ one has 
\begin{eqnarray}
&&e^{q(t-s)}>\frac{q(e^{-\mu(t-s)}-e^{-\mu t})+\mu(e^{q (t-s)}-e^{-\mu t})}{(\mu+q)(1-e^{-\mu t})}
\label{eq24}\\
&\Leftrightarrow&qe^{q(t-s)}+(\mu+q)e^{-\mu t}-(\mu+q)e^{q(t-s)-\mu t}-qe^{-\mu(t-s)}> 0
\nn\\
&\Leftrightarrow& e^{q(t-s)}f(t)>0,\;\;{\rm where}
\nn\\
f(t)&:=&q+(\mu+q)e^{-q(t-s)-\mu t}-(\mu+q)e^{-\mu t}-qe^{-(q+\mu)(t-s)}.
\nn
\end{eqnarray}
Then 
\begin{eqnarray}
f(s)&=&0,
\nn\\
f'(t)&=&(\mu+q)e^{-\mu t}[-(\mu+q)e^{-q(t-s)}+\mu+qe^{-q(t-s)+\mu s}]
\nn\\
&=&(\mu+q)e^{-\mu t}[qe^{-q(t-s)}(e^{\mu s}-1)+\mu(1-e^{-q(t- s)})]>0.
\nn
\end{eqnarray}
Hence $f(t)>0$ for all $t>s$ and (\ref{eq24}) holds. Setting $s=\t$ and
$s=0$ shows the respective statements. 
\qed

%
\section{The stem cell model formulated as a SD-DDE}\label{s1}
In this section, in regard to Section \ref{s3} and the DDE (\ref{eq18}), we keep the assumptions on $q$, $\mu$ and $\cD$, but specify $j$  and $\t$, such that the DDE (\ref{eq18})
becomes the SD-DDE (\ref{eq11}--\ref{eq14}) that describes the stem cell 
dynamics. Then we apply the previous results to analyze this SD-DDE. 
%
\subsection{Assumptions and main results}\label{ss1}
Suppose that the function $g$ satisfies the following property, which we denote 
by (G):
There exist $x_1, x_2, b, K,\e \in\R$, such that $x_1<x_2$, $0<\e<K$ and $b>0$,
and $g:\oB_b(x_2)\times \R_+\->\R$
\begin{itemize}
\item[($G_1$)] is locally Lipschitz in the second argument, uniformly with respect to the first,
\item[($G_2$)] is partially differentiable with respect to the first argument with $D_1g$
Lipschitz and 
\[
\sup_{(y,z)\in \oB_b(x_2)\times \R_+}|D_1g(y,z)|<\frac{K}{b},
\]
\item[($G_3$)] satisfies $\e\le g(y,z)\le K$ on $\oB_b(x_2)\times \R_+$ and $x_2-x_1\in(0,\frac{b}{K}\e)$.
\end{itemize}
Note that $(G_3)$ implies that $x_1\in\oB_b(x_2)$. We now define 
$h:=\frac{b}{K}$. The following result  is an application of the Picard-Lindel\"of 
theorem.
\begin{lemma}
Let $\psi\in C([-h,0],\R_+)$. Then there exists a unique solution $y=y(\cdot,\psi)$  on $[0,h]$ of
(\ref{eq13})
with $y([0,h],\psi)\subset\oB_b(x_2)$. Moreover, there exists a unique 
\[
\t=\t(\psi)\in[\frac{x_2-x_1}{K},\frac{x_2-x_1}{\e}]\subset(0,h)
\]
solving (\ref{eq14}). 
\end{lemma}
\Proof Define $f_\psi:[0,h]\times\oB_b(x_2)\->\R;\;f_\psi(s,y):=-g(y,\psi(-s))$ and with $f_\psi$
a non-autonomous ODE $y'(s)=f_\psi(s,y(s))$. Then (G)  guarantees directly that $f_\psi$
satisfies the conditions of the Picard-Lindel\"of Theorem, e.g. \cite[Theorem II.1.1]{Hartman}, 
which guarantees that there exists a unique solution $y$ on $[0,h]$, since we defined $h:=\frac{b}{K}$. The remaining statements can be shown by integrating the ODE and using (G3). 
\qed

Accordingly, with $\uta:=(x_2-x_1)/K$ we can now define a functional 
\[
\t:C([-h,0],\R_+)\->[\uta,h)
\] 
to describe the state-dependence of the delay. Moreover, we
suppose that $d:\oB_b(x_2)\times \R_+\->\R$ is bounded and Lipschitz and that 
$\g:\R_+\->\R_+$ is bounded and locally Lipschitz and define
\begin{eqnarray}
j(\vi,\psi):=
\frac{\g(\psi(-\t(\psi)))}{g(x_1,\psi(-\t(\psi)))}g(x_2,\psi(0))\vi(-\t(\psi))e^{\int_0^{\t(\psi)}
[d-D_1g](y(s,\psi),\psi(-s))ds}.
\nn\\
\label{eq17}
\end{eqnarray}
Then clearly the DDE (\ref{eq18}) becomes the SD-DDE (\ref{eq11}--\ref{eq14}) and (\ref{eq20}) holds 
with
\begin{eqnarray}
k_j:=\frac{K}{\e}
e^{(\frac{K}{b}+\sup_{(y,z)\in\oB_b(x_2)\times\R_+}|d(y,z)|)h}\sup_{z\in\R_+}\g(z)<\infty. 
\nn
\end{eqnarray}
The following result will be proven in the next subsection.
\begin{theorem}\label{theo9}
For any $\phi=(\vi,\psi)\in V_\cD$, under the conditions given in this subsection, the SD-DDE (\ref{eq11}--\ref{eq14}) has a unique
solution $x^\phi=(w,v)$ on $\R_+$ through $\phi$. The solutions
define a continuous semiflow in the 
sense of Theorem \ref{theo8} and with $f_\t$ as in Theorem \ref{theo10} (c) satisfy the invariance properties Theorem \ref{theo10} (e-f). 
\end{theorem}
%
%
%
\subsection{Proofs}\label{ss2}
We can apply Theorem \ref{theo10} to obtain the statement of Theorem \ref{theo9} if we show that $j$ is almost locally Lipschitz. To show this, it is useful to introduce a notation that summarizes
model ingredients with the same type of delay:
Let first $\b:\R_+\->\R$, $r:C([-h,0],\R_+)\->[0,h]$ and $\cG:C([-h,0],\R_+)\->\R$ be arbitrary maps.  
As a tool to prove several results that follow we define the evaluation operator 
\begin{eqnarray}
C([-h,0],\R_+)\times[-h,0]\->\R;\;ev(\vi,s):=\vi(s).
\label{eq15}
\end{eqnarray}
 Trivially, $ev$ inherits continuity from the functions in its domain. We will show that $j$ is a special
 case of the functional defined in the following lemma. 
\begin{lemma}\label{lem14}
Suppose that $\b$ is locally Lipschitz and that $r$ and $\cG$ are almost locally Lipschitz, then 
the functional $\cD\->\R_+;$
\begin{eqnarray}
(\vi,\psi)\longmapsto\b(\psi(-r(\psi)))\vi(-r(\psi))\cG(\psi)
\label{eq16}
\end{eqnarray}
 is almost locally Lipschitz. 
\end{lemma}
\Proof By the discussed sum - and product rules and by other rules, which are straightforward, it suffices to show that
the two maps $\psi\mapsto\b(\psi(-r(\psi)))$ and $(\vi,\psi)\mapsto\vi(-r(\psi))$ are almost locally Lipschitz. Now note that the first map can be decomposed as
\[
\psi\mapsto(\psi,-r(\psi))\xmapsto{ev}\psi(-r(\psi))\xmapsto{\b}\b(\psi(-r(\psi))). 
\]
Hence it is continuous as a composition by continuity of $r$, $ev$ and $\b$. Similarly the
second map can be written as 
\[
(\vi,\psi)\mapsto(\vi,-r(\psi))\xmapsto{ev}\vi(-r(\psi))
\]
and continuity can be concluded. Next, let $\psi_0\in C([-h,0],\R_+)$, $R>0$. Choose $\de>0$, 
$k$ such that $r$ is $k$-Lipschitz on $V(\psi_0;\de,R)$. Now note that for $\psi,\chi\in V(\psi_0;\de,R)$
\begin{eqnarray}
&&|\psi(-r(\psi)-\chi(-r(\chi)))|
\nn\\
&\le&|\psi(-r(\psi)-\psi(-r(\chi)))|+|\psi(-r(\chi)-\chi(-r(\chi)))|
\nn\\
&\le&R|r(\psi)-r(\chi)|+\|\psi-\chi\|\le(Rk+1)\|\psi-\chi\|.
\nn
\end{eqnarray}
Hence, $\psi\mapsto\psi(-r(\psi))$ is almost locally Lipschitz. Since $\b$ is locally Lipschitz,
$\psi\mapsto\b(\psi(-r(\psi)))$ is almost local Lipschitz by the discussed composition rule. The stated
Lipschitz property of the second map follows similarly. 
\qed

A Gronwall-Lemma type estimate and use of ($G_1$) and ($G_2$) lead to the following result. 
\begin{lemma}\label{lem8}
The map $Y: C([-h,0],\R_+)\->C([0,h],\oB_b(x_2));\;Y(\psi)(t):=y(t,\psi)$ is locally Lipschitz. 
\end{lemma}
\Proof Let $\psi_0,\psi, \ops \in C([-h,0],\R_+)$. 
One has
\begin{eqnarray}
&&|g(y(s,\psi),\psi(-s))-g(y(s,\ops),\ops(-s))|
\nn\\
&\le&|g(y(s,\psi),\psi(-s))-g(y(s,\ops),\psi(-s))|
\nn\\
&&+|g(y(s,\ops),\psi(-s))-g(y(s,\ops),\ops(-s))|
=:(I)+(II)
\nn
\end{eqnarray}
in obvious notation. By ($G_2$) and the mean value theorem one has
\begin{eqnarray}
(I)&\le&L_1|y(s,\psi)-y(s,\ops)|\le L_1\|y(\cdot,\psi)-y(\cdot,\ops)\|= L_1\|Y(\psi)-Y(\ops)\|,
\nn
\end{eqnarray}
where $L_1:=\sup_{(y,z)\in\oB_b(x_2)\times \R_+}|D_1g(y,z)|$. 
By ($G_1$), one has $(II)\le L_2\|\psi-\ops\|$ for some $L_2\ge 0$ and $\psi$ and $\ops$ in a neighborhood of $\psi_0$. Now combine
\begin{eqnarray}
|y(t,\psi)-y(t,\ops)|\le \int_0^t|g(y(s,\psi),\psi(-s))-g(y(s,\ops),\ops(-s))ds|
\nn
\end{eqnarray}
with the previous estimates and $L_1<\frac{1}{h}$, which follows from ($G_2$), to complete the proof. 
\qed

We can use this result to deduce
\begin{lemma}\label{lem9}
The map $C([-h,0],\R_+)\->[0,h];\;\psi\mapsto\t(\psi)$ is locally Lipschitz. 
\end{lemma}
\Proof Let $\ops, \psi\in C([-h,0],\R_+)$. By definition of $\t(\psi)$ and $\t(\ops)$ one has
\begin{eqnarray}
y(\t(\psi),\psi)=y(\t(\ops),\ops)\;\;(=x_1).
\nn
\end{eqnarray}
Hence, 
\begin{eqnarray}
|y(\t(\psi),\psi)-y(\t(\psi),\ops)|=|y(\t(\psi),\ops)-y(\t(\ops),\ops)|.
\nn
\end{eqnarray}
The left hand side is dominated by $\|Y(\psi)-Y(\ops)\|$. There exists some $t\in[0,h]$, such
that the right hand side equals
\begin{eqnarray}
&&|D_1y(t,\ops)||\t(\psi)-\t(\ops)|
\nn\\
&=&|g(y(t,\ops),\ops(-t))||\t(\psi)-\t(\ops)|\ge\e|\t(\psi)-\t(\ops)|
\nn
\end{eqnarray}
by $(G_3)$. Thus $|\t(\psi)-\t(\ops)|\le\frac{1}{\e}|Y(\psi)-Y(\ops)|$ and the proof
 is completed using Lipschitzianity of $Y$. 
 \qed

 \begin{lemma}\label{lem10}
Let $G:C([-h,0],\R_+)\times C([0,h],\oB_b(x_2)) \->C([0,h],\R)$ be an arbitrary locally Lipschitz operator with
\[
\sup_{(\psi,z)}lip\;G(\psi,z)<\infty.
\]
Define $\cG:C([-h,0],\R_+)\->\R;\;\cG(\psi):=g(x_2,\psi(0))e^{G(\psi,Y(\psi))(\t(\psi))}$. Then $\cG$ 
is locally Lipschitz. 
 \end{lemma}
 \Proof Choose $\vi_0\in C([-h,0],\R_+)$ and $R:=\sup_{(\psi,z)}lip\;G(\psi,z)$. Choose
$k$ and $\de$ such that $G$ is $k$-Lipschitz on $\oB_\de((\vi_0,Y(\vi_0)))$ and $Y$ and $\t$
are $k$-Lipschitz on $\oB_\de(\vi_0)$. Choose $\e\le\de$ such that
\[
|G(\psi,Y(\psi))-G(\vi_0,Y(\vi_0))|\le\de,\;\;{\rm if}\;\psi\in\oB_\e(\vi_0). 
\]
Let $\vi,\psi\in\oB_\e(\vi_0)$. Now note that $ev$ is $\max\{\overline R,1\}$-Lipschitz on $V(\ovi;\overline\de,\overline R)\times[-h,0]$ for any $\ovi$, $\overline\de$, $\overline R$. Hence, 
for $\psi, \chi\in\oB_\e(\vi_0)$ 
\begin{eqnarray}
&&|ev(G(\psi,Y(\psi)),\t(\psi))-ev(G(\chi,Y(\chi)),\t(\chi))|
\nn\\
&\le&\max\{R,1\}\{|G(\psi,Y(\psi))-G(\chi,Y(\chi))|+|\t(\psi)-\t(\chi)|\}
\nn\\
&\le&\max\{R,1\}\{k[\|\psi-\chi\|+\|Y(\psi)-Y(\chi)\|]+|\t(\psi)-\t(\chi)|\}
\nn\\
&\le&\max\{R,1\}\max\{k,k^2\}\|\psi-\chi\|.
\nn
\end{eqnarray}
Thus $\psi\mapsto ev(G(\psi,Y(\psi)),\t(\psi))$ is locally Lipschitz. This implies that $\cG$
is locally Lipschitz. 
 \qed
 
 \begin{lemma}\label{lem13}
Let $J\subset\R$, $k:J\times \R_+\->\R$ be an arbitrary Lipschitz and bounded map. Define 
\begin{eqnarray}
&&G:C([-h,0],\R_+)\times C([0,h],J)
\->C([0,h],\R);
\nn\\
&&G(\psi,z)(t):=\int_0^tk(z(s),\psi(-s))ds.
\nn
\end{eqnarray}
Then $G$ is Lipschitz and
\[
\sup_{(\psi,z)}lip\;G(\psi,z)<\infty.
\] 
 \end{lemma}
 \Proof The first result follows from the estimates
\begin{eqnarray}
&&|G(\psi,z)(t)-G(\ops,\oz)(t)|\le\int_0^t|k(z(s),\psi(-s))-k(\oz(s),\ops(-s))|ds, 
\nn\\
&&|k(z(s),\psi(-s))-k(\oz(s),\ops(-s))|\le L[\|z-\oz\|+\|\psi-\ops\|],
\nn
\end{eqnarray}
for some $L\ge 0$.
Boundedness of $k$ implies the second statement. 
\qed

\begin{remark}
In Lemma \ref{lem10} and below we merely need local Lipschitzianity of $G$. We presented a 
sketch of the rather straightforward proof of the previous lemma to also hint that
mere local Lipschitzianity of $k$ would not yield local Lipschitzianity of G, however. 
The point is that continuous functions being close in a point obviously in general does not make them close in the sup-norm. See \cite{Appell} for more details on smoothness
properties of related Nemytskii-operators. 
\end{remark}
 We now combine our results and prove
\begin{prop}\label{prop1}
The functional $j$ as defined in (\ref{eq17}) is almost locally Lipschitz. 
\end{prop}
\Proof By Lemma \ref{lem14} it is sufficient to show
that $\frac{\g(\cdot)}{g(x_1,\cdot)}$ is locally Lipschitz and that $\t$ and $\psi\mapsto g(x_2,\psi(0))\exp\{\int_0^{\t(\psi)}[d-D_1g](y(s,\psi),\psi(-s))ds\}$ are almost locally
Lipschitz. Local Lipschitzianity of the first map follows directly from local Lipschitzianity of
$\g$, ($G_1$) and ($G_3$). 
(Almost) local Lipschitzianity of $\t$ is shown in Lemma \ref{lem9}. (Almost) local Lipschitzianity
of the third map follows by Lemma \ref{lem10}, provided we show local Lipschitzianity
of $G$, defined as $G(\psi,z)(t):=\int_0^t[d-D_1g](z(s),\psi(-s))ds$ and that for this $G$ one
has $\sup_{(\psi,z)}lip\;G(\psi,z)<\infty$. The latter follow
by Lemma \ref{lem13} from boundedness and Lipschitzianity of $k:=d-D_1g$ with $J:=\oB_b(x_2)$. Thus, $j$ is almost locally
Lipschitz. %
\qed

%
\section{Examples of model ingredients}\label{ss5}
In the previous section we have elaborated conditions on the model ingredients specified as functions $q$,  $\g$, $g$ and $d$ and the nonnegative parameter $\mu$. The exact nature of the cellular and sub-cellular processes related to
these ingredients is subject to current research  \cite{Vivanco}. In \cite{Getto1} a combination of available knowledge with mathematical considerations led to the specification
\begin{eqnarray}
&&q(z):=[2s_w(z)-1]d_w(z)-\mu_w,\;\;\g(z):=2[1-s_w(z)]d_w(z),\;\;{\rm where}
\nn\\
&& s_w(z):=\frac{a_w}{1+k_a z},\;\;d_w(z):=\frac{p_w}{1+k_pz}
\nn
\end{eqnarray}
with $a_w\in[0,1]$ and $p_w$, $\mu_w$, $k_a$ and $k_p$ nonnegative parameters. It
is obvious that for these examples $q$ and $\g$ are Lipschitz, in particular locally
Lipschitz. The function $d$ considered is of the form
\[
d(y,z)=\frac{\a(y)}{1+k_dz}-\mu_u(y)
\]
for a nonnegative parameter $k_d$ and nonnegative functions $\a$ and $\mu_u$. Note that we
here assumed the $y$-component of the domain to be compact ($\oB_b(x_2)$). Hence, if $\a$ and $\mu_u$ are Lipschitz, then $d$ is Lipschitz and bounded. 

In \cite{Doumic} based on \cite{Stiehl} the authors consider $g$ of the shape
\begin{eqnarray}
g(y,z)=2[1-\frac{a(y)}{1+k_gz}]p(y)
\label{eq26}
\end{eqnarray}
for nonnegative $k_g$, $a$ and $p$. 
Further specifications are considered, which lead to $y$- and $z$-independent $g$  respectively. 
We here suppose that $a$ and $p$ are differentiable and that $a'$ and $p'$ are Lipschitz. If we slightly modify (\ref{eq26}) such that $g(y,z)\ge\e$ on $\oB_b(x_2)\times\R_+$, and choose the
constants in (G) appropriately, we can guarantee that $g$ satisfies (G). 

Note that, though our assumption that $g$ is bounded away from zero
has a mathematical motivation, a nonzero maturation rate also has biological consistency. An example of a $g$ that is decreasing in $z$ could be
\begin{eqnarray}
g(y,z):=\e+e^{-z}\g_g(y) 
\nn
\end{eqnarray}
with $\g_g$ differentiable and $\g_g'$ Lipschitz.

 A choice $g(y,z)\equiv 1$ also fulfills the requirements and with this choice $y$
 could be interpreted as the age of a progenitor cell. 
%
\section{Discussion and outlook}\label{s6}

Note that in \cite[Theorem 6.8]{Nishi} a large class of SD-DDE is analyzed. An alternative approach to proving well-posedness for (1.1-1.4) could be, to investigate whether the cited result can be modified to include distributed delays and whether the there required smoothness conditions can be guaranteed. Possibly also with that approach the implementation
of retractions could be useful. For results on differentiability of solutions with respect to parameters and initial data, which are related to our results on continuous dependence on
initial values, we refer to the work of Hartung, e.g. \cite{Hartung1, Hartung2}. 

For the specifications in Section \ref{ss5} and under some additional assumptions, see \cite{Getto1}, the here analyzed model (\ref{eq11}--\ref{eq14}) has a unique positive equilibrium emerging from the trivial equilibrium in a transcritical bifurcation:  the rate $q$ can be assumed to be
decreasing to a negative value, hence the bifurcation parameter should guarantee that $q(0)>0$. 

In a manuscript in preparation Ph.G. and G.R.
are using the theory of  \cite{Smith} to show that the trivial equilibrium is globally 
asymptotically stable in absence of the positive equilibrium, whereas in its 
presence, there is uniform strong population persistence. The latter can be 
concluded, essentially, if the system is dissipative. 
In the manuscript, Ph.G. and G.R.
 encounter a situation in which there either is dissipativity or $\cC$-convergence of the solution to a constant, where the constant depends on the initial condition. A priori it is not clear how dissipativity can be concluded from
the second case. By Corollary \ref{corol1}, however, it can be concluded that the constant is an equilibrium solution and, as the equilibrium is unique, this implies dissipativity.
 Note also that Corollary \ref{corol1} follows from continuous dependence of the solution on the initial value in the $\cC$-topology. 
 Using continuous dependence of the solution on the initial value in the $\cC^1$-topology, as established in \cite{Getto}, one could possibly prove similarly that the limit is an equilibrium, if the convergence of the solution to the constant is in $\cC^1$. 
In the manuscript in preparation, the authors, however, are not able to show this convergence in  $\cC^1$. 
Hence, a $\cC^1$-variant of Corollary \ref{corol1} would not be applicable in that manuscript.  In summary
the present Corollary \ref{corol1} can be expected to be a necessary and  sufficient tool to show dissipativity and uniform strong persistence for 
(\ref{eq11}--\ref{eq14}).

In \cite{Getto}, the derivative of the semiflow defined on the solution manifold is
computed, such that a linearization is at hand. General theorems of linearized stability, applicable to our system, are shown in \cite{Walther}  (stability) and \cite{Stumpf} (instability). 
By the analysis of the characteristic equation derived from this linerization
in a manuscript in preparation by Mats Gyllenberg, Yukihiko Nakata, Francesca Scarabel and Ph. G., the positive equilibrium is stable upon emergence in the neighborhood of a transcritical bifurcation point and destabilizes by a pair of eigenvalues crossing into the right half plane. Based on this analysis and on unpublished numerical simulations with DDE-biftool \cite{Sieber} (by Jan Sieber) and pseudo-spectral methods \cite{Breda} (by F. Scarabel)
there is evidence for a Hopf bifurcation and the emergence of a limit cycle. 
This motivates the idea of a future analysis of Hopf bifurcations and periodic solutions.

We refer to \cite{Huwu} for Hopf bifurcation analysis for related equations. 
To establish periodicity for a general class of equations, in \cite{MN1} the authors include
the assumption that the initial function should be at equilibrium value at time zero. 
If for our model this assumption is included, one can investigate convex and compact sets that are invariant under the original untransformed system (\ref{eq11}-\ref{eq14}), i.e., sets that are invariant for both components of the state. 
Motivated by the fact that periodicity for infinite times often can be concluded from behavior in a finite time interval, we also have some hope that the here
established invariance for finite time may be sufficient. 
\bigskip

\noindent
{\bf Acknowledgements:} The manuscript was inspired by discussions with
Tibor Krisztin during a postdoctoral stay of Ph.G. at the University
of Szeged. Ph.G. thanks Stefan Siegmund und Reinhard Stahn at Technische Universit\"at Dresden for help with the manuscript. 
%



\end{document}